\documentclass[11pt,titlepage,a4paper]{article}
\usepackage{bozhomac}   
\usepackage{bozhlogo}   
\usepackage{varioref}	
\usepackage{cite}   	
\usepackage{tabularx}   
\usepackage{pb-diagram}	
\usepackage{fancybox}   

\title{\bfseries    \vspace*{-1.678902345in}
{\huge Links between connections, parallel
\\[1.22ex] transports and transports along paths
}
}

\vspace{1.7ex}

\author{
Bozhidar Z.\ Iliev
\thanks{Laboratory of Mathematical Modeling in Physics,
Institute for Nuclear Research and \mbox{Nuclear} Energy,
Bulgarian Academy of Sciences,
Boul.\ Tzarigradsko chauss\'ee~72, 1784 Sofia, Bulgaria}
\thanks{E-mail address: bozho@inrne.bas.bg}
\thanks{URL: http://theo.inrne.bas.bg/$\sim$bozho/}
}

\date{
 \vspace{2.27ex}\ShortTitle{Connections and (parallel) transports along paths}
								\\[0.27ex]
 \vspace{3.27ex}
\small
    \begin{tabular}{r@{$\colon\to~$}l}
 \vspace{0.09ex} Last update   & February 29, 2008	\\[0.09ex]
    \end{tabular} \\[1.27ex]
\normalsize
\vspace{0.27ex}
\textsl{\bfseries
Report presented at the International Workshop\\
``Advanced Geometric Methods in Physics''\\
Florence, Italy, 14 -- 18 April, 2005}		\\[2ex] %
\small
    \begin{tabular}{r@{$\colon~$}l}
\normalsize\sffamily\bfseries
 \vspace{0.27ex} http://www.arXiv.org e-Print archive No. &
\normalsize\sffamily\bfseries
math.DG/0504010				\\[1.27ex]
    \end{tabular} \\[-0.27ex]
\normalsize
 \vspace{4.27ex}{\Huge \BOZHO}  \\[4.27ex]
\vspace{0.27ex}\Subject{Differential geometry}    \\[2.27ex]
    \begin{tabular}{r@{\hspace{0.512em}}|@{\hspace{0.512em}}l}
 \vspace{0.27ex}
&
 \vspace{0.27ex}
    \end{tabular} \\[1.27ex]
 \vspace{0.27ex}    \\[0.27ex]
}

\begin{document}        

\renewcommand{\thepage}{\roman{page}}

\renewcommand{\thefootnote}{\fnsymbol{footnote}} 
\maketitle              
\renewcommand{\thefootnote}{\alph{footnote}}   

\tableofcontents        


\begin{abstract}

The axiomatic approach to parallel transport theory is partially discussed.
Bijective correspondences between the sets of connections, (axiomatically
defined) parallel transports and transports along paths satisfying some
additional conditions are constructed. In particular, the equivalence
between the concepts ``connection'' and ``(axiomatically defined) parallel
transport'' is established. Some results which are specific for topological and
vector bundles are presented.

\end{abstract}

\renewcommand{\thepage}{\arabic{page}}


\section {Introduction}
\label{Introduction}

	 From different view\ndash points, the connection theory can be found
in many works, like%
~\cite{Kobayashi-1957,K&N-1,Sachs&Wu,Nash&Sen,Warner,Bishop&Crittenden,
Yano&Kon,Steenrod,Sulanke&Wintgen,Bruhat,Husemoller,Mishchenko,R_Hermann-1,
Greub&et_al.-1,Atiyah,Dandoloff&Zakrzewski,Tamura,Hicks,Sternberg,Rahula,
Mangiarotti&Sardanashvily}.
As pointed in these and many other references, the concept of a parallel
transport is defined on the base of the one of a connection. The opposite
approach, \ie the definition of a connection on the ground of an axiomatically
defined one of a parallel transport, is considered
in~\cite{Lumiste-1964,Lumiste-1966,Teleman,Dombrowski,Lumiste-1971,
Mathenedia-4,Durhuus&Leinaas, Khudaverdian&Schwarz,Ostianu&et_al.,Nikolov,
Poor}. The major classical results on axiomatization of parallel transport
theory are presented in the table below

\renewcommand{\arraystretch}{1.32}
    \begin{table}[ht!]  \label{HistoricalTablePT}%
\index{parallel transport!history of}
    \begin{tabularx}{\textwidth}{@{}rlX@{}}
Year    & Person    & Result and original reference
\\ \hline
1917    & T. Levi-Civita%
            & Definition of a parallel transport of a vector in Riemannian
geometry.~\cite{Levi-Civita/1917}
\\
1949    & Willi Rinow%
            & An axiomatic definition of parallel transport in tangent bundle
is introduced in unpublished lectures at Humboldt
university. (See~\cite{Dombrowski} and~\cite[p.~46]{Poor}.)
\\
1964    & \"{U}.G. Lumiste%
            & Definition of a connection in principal bundle (with homogeneous
fibres) as a parallel transport along canonical paths $\alpha\colon[0,1]\to M$
in its base $M$. The parallel transport is defined as a mapping from the fibre
over $\alpha(0)$ into the one over $\alpha(1)$ satisfying some
axioms.~\cite[sec.~2.2]{Lumiste-1964}
\\
1964    & C. Teleman%
            & Definition of a connection in topological bundle as a parallel
transport along canonical paths $\alpha\colon[0,1]\to M$ in its base $M$. The
parallel transport is defined as a lifting  of these paths through a point in
the fibre over their initial points
$\alpha(0)$.~\cite[chapter~IV, sec.~B.3]{Teleman}
\\
1968
& P. Dombrowski%
            & Definition of a linear connection in vector bundle as a parallel
transport along paths $\beta\colon[a,b]\to M$, with $a,b\in\field[R]$ and
$a\le b$, in its base $M$. The parallel transport is defined as a mapping from
the fibre over $\alpha(a)$ into the one over $\alpha(b)$ satisfying certain
axioms. The theory of covariant derivatives is constructed on that
base.~\cite[\S~1]{Dombrowski}
\\
1981    & Walter Poor%
            & A detailed axiomatic definition of a parallel transport in vector
bundles. The whole theory of linear connections in such bundles is deduced on
that ground.\cite{Poor}
\\\hline
    \end{tabularx}
    \end{table}

	The purpose of the present investigation is to be revealed different
relations between the concepts ``connection on a bundle'' and
``(parallel) transport (along paths) in a bundle''.
These concepts will be shown to be equivalent in a sense that there exist
bijective mappings between the sets of connections, (axiomatically defined)
parallel transports and transports along paths satisfying some additional
conditions. In this way, two different, but equivalent, systems of axioms
defining the concept ``parallel transport'' will be established.

	The paper is organized as follows.

	Section~\ref{Sect3} recalls the definition of a connection on a
bundle  as a distribution on its bundle space which is complimentary to the
vertical distribution on it. The notions of parallel transport generated by
connection  and of (axiomatically defined) parallel transport are introduced.
	In section~\ref{Sect4} is defined the notion of a transport along
paths, which is one of the possible axiomatic approaches to the concept of a
``parallel transport''.
	Section~\ref{Sect5} shows how a transport along paths generates a
connection and presents some properties of such a connection.
	The main result of section~\ref{Sect6} is the construction of
bijective mapping between the sets of transports along paths, satisfying some
additional conditions, and parallel transports.
	Section~\ref{Sect7} deals with different relations between
connections, parallel transports, and transports along paths. The
equivalence between the concepts connection, transport along paths,
satisfying some additional conditions, and parallel transport is
established.
	Section~\ref{Sect8} summarizes the basic results of the present work
in a series of theorems.
	In section~\ref{Subsect8.5} some of the preceeding results are
generalized in a case of general topological bundles.
	In section~\ref{Subsect8.6} are proved some results which are specific
for vector bundles. In particular, a conditions ensuring equivalence between
linear transport along paths and a linear connection in these bundles is
presented.
	Section~\ref{Conclusion} closes the paper with some concluding
remarks.

 \vspace{1.25ex}

        Let us now present some introductory material, like notation etc.,
that will be needed for our exposition. The reader is referred for details to
standard books on differential geometry, like~\cite{K&N,Warner,Poor}.

        A differentiable finite-dimensional manifold over a field $\field$
will be denoted typically by $M$. Here $\field$ stands for the field
$\field[R]$ of real or the field $\field[C]$ of complex numbers,
$\field=\field[R],\field[C]$. The manifolds we consider are supposed to be
smooth of class $C^1$.~%
\footnote{~%
Some of our definitions or/and results are valid also for $C^0$
manifolds, but we do not want to overload the material with continuous
counting of the required degree of differentiability of the manifolds
involved. Some parts of the text admit generalizations on more general
spaces, like the topological ones, but this is out of the subject of the
present work.%
}

The set of vector fields, realized as first order differential operators,
over $M$ will be denoted by $\mathcal{X}(M)$. The space tangent to $M$ at
$p\in M$ is $T_p(M)$  and $(T(M),\pi_T,M)$ stand for the tangent bundle
over $M$. The value of $X\in\mathcal{X}(M)$ at $p\in M$ is $X_p\in T_p(M)$.

	By $J\subseteq\field[R]$ will be denoted an arbitrary real interval
that can be open or closed at one or both its ends.
	The notation $[\sigma,\tau]$ will be used for an arbitrary closed
interval with ends $\sigma,\tau\in\field[R]$, with $\sigma\le \tau$.
	The notation $\gamma\colon J\to M$ represents an arbitrary path in $M$.
	For a $C^1$ path $\gamma\colon J\to M$, the vector tangent to $\gamma$
at $s\in J $ will be denoted by
\(
\dot\gamma(s)
:=\frac{\od} {\od t}\Big|_{t=s}(\gamma(t)) \in T_{\gamma(s)}(M) .
\)
If $s_0\in J$ is an end point of $J$ and $J$ is closed at $s_0$, the
derivative in the definition of $\dot{\gamma}(s_0)$ is regarded as a
one\ndash sided derivative at $s_0$.

	In this work, we shall need the notion of an inverse path and of a
product of paths. There are not `natural' definitions of these concepts, but
this is not important for us as the (parallel) transports we shall consider
bellow are parametrization invariant in some sense, like~\eref{8.6} below. For
that reason, the concepts mentioned will be defined only for \emph{canonical
paths} $[0,1]\to M$, whose domain is the real interval $[0,1]:=\{r\in\field[R]
: 0\le r \le 1 \}$.
    The path \emph{inverse} to $\gamma\colon[0,1]\to M$ is
$\gamma_{\_}:=\gamma\circ\tau_{\_}\colon [0,1]\to M$, with $\circ$ being the
sign of composition of mappings and $\tau_{\_}\colon[0,1]\to[0,1]$ being
given by $\tau_{\_}(t):=1-t$ for $t\in[0,1]$.
    If $\gamma_1,\gamma_2\colon[0,1]\to M$ and $\gamma_1(1)=\gamma_2(0)$,
the \emph{product} $\gamma_1\gamma_2$ of $\gamma_1$ and $\gamma_2$ is a
canonical path $\gamma_1\gamma_2\colon[0,1]\to M$ such that
 $(\gamma_1\gamma_2)(t):=\gamma_1(2t)$ for $t\in[0,1/2]$ and
 $(\gamma_1\gamma_2)(t):=\gamma_2(2t-1)$ for $t\in[1/2,1]$.
For more details on this item, see~\cite{Sze-Tsen,Viro&Fuks}.

        If $k\in \field[N]$ and $k\le\dim M$, a $k$-dimensional
\emph{distribution} $\Delta$ on $M$ is defined as a mapping
$\Delta\colon p\mapsto \Delta_p$ assigning to each $p\in M$ a $k$\ndash
dimensional subspace $\Delta_p$ of the tangent space $T_p(M)$ of $M$ at $p$,
$\Delta_p\subseteq T_p(M)$.
	We say that a vector field $X\in\mathcal{X}(M)$ is in $\Delta$ and
write $X\in\Delta$, if $X_p\in\Delta_p$ for all $p\in M$.

	By $(E,\pi,M$), we shall denote a bundle with bundle space $E$,
projection $\pi\colon E\to M$, and base space $M$.
	Suppose that the spaces $M$ and $E$ are
manifolds of finite dimensions $n\in\field[N]$ and $n+r$, for some
$r\in\field[N]$, respectively; so the dimension of the fibre $\pi^{-1}(x)$,
with $x\in M$, \ie the fibre dimension of $(E,\pi,M)$, is $r$.
Besides, let these manifolds be $C^1$ differentiable, if the opposite is not
stated explicitly.~%
\footnote{~%
Some of our considerations are valid also for $C^0$ manifolds. By assuming
$C^1$ differentiability, we skip the problem of counting the required
differentiability class of the whole material that follows. Sometimes, the
$C^1$ differentiability is required explicitly, which is a hint that a
statement or definition is not valid otherwise.%
}


\section
[Connections and the generated by them parallel transports (review)]
{Connections and the generated by them parallel transports\\ (review)}
\label{Sect3}

	From a number of equivalent definitions of a connection on
differentiable manifold~\cite[sections~2.1
and~2.2]{Mangiarotti&Sardanashvily}, we shall use the following one.

    \begin{Defn}    \label{Defn3.1}
A \emph{connection on a bundle} $(E,\pi,M)$ is an $n=\dim M$
dimensional distribution $\Delta^h$ on $E$ such that, for each $p\in E$
and the \emph{vertical distribution} $\Delta^v$ defined by
    \begin{equation}    \label{3.9-2}
\Delta^v\colon p \mapsto \Delta^v_p
:= T_{\imath(p)}\bigl( \pi^{-1}(\pi(p)) \bigr)
\cong T_{p}\bigl( \pi^{-1}(\pi(p)) \bigr),
    \end{equation}
with $\imath\colon\pi^{-1}(\pi(p))\to E$ being the inclusion mapping, is
fulfilled
    \begin{equation}    \label{3.9-3}
\Delta^v_p\oplus \Delta^h_p = T_p(E) ,
    \end{equation}
where
\(
\Delta^h\colon p \mapsto \Delta^h_p \subseteq T_{p}(E)
\)
and $\oplus$ is the direct sum sign. The distribution $\Delta^h$ is
called \emph{horizontal} and symbolically we write
$\Delta^v\oplus\Delta^h=T(E)$.
    \end{Defn}

    A \emph{vector} at a point $p\in E$ (resp. a \emph{vector field} on $E$)
is said to be \emph{vertical} or \emph{horizontal} if it (resp.\ its value at
$p$) belongs to $\Delta^v_p$ or $\Delta^h_p$, respectively, for the given
(resp.\ any) point $p$.
	A vector $Y_p\in T_p(E)$ (resp.\ vector field $Y\in\mathcal{X}(E)$) is
called a \emph{horizontal lift of a vector} $X_{\pi(p)}\in T_{\pi(p)}(M)$
(resp.\ vector \emph{field} $X\in\mathcal{X}(M)$ on $M=\pi(E)$) if
$\pi_*(Y_p)=X_{\pi(p)}$ for the given (resp.\ any) point $p\in E$. Since
$\pi_*|_{\Delta_p^h}\colon \Delta_p^h\to T_{\pi(p)}(M)$ is a vector space
isomorphism for all $p\in E$~\cite[sec.~1.24]{Poor}, any vector in
$T_{\pi(p)}(M)$ (resp.\ vector field in $\mathcal{X}(M)$) has a unique
horizontal lift in $T_p(E)$ (resp.\ $\mathcal{X}(E)$).

	As a result of~\eref{3.9-3}, any vector $Y_p\in T_p(E)$ (resp.\ vector
field $Y\in\mathcal{X}(E)$) admits a unique representation
 $Y_p=Y_p^v\oplus Y_p^h$ (resp.\ $Y=Y^v\oplus Y^h$) with
 $Y_p^v\in\Delta_p^v$ and $Y_p^h\in\Delta_p^h$
(resp.\ $Y^v\in\Delta^v$ and $Y^h\in\Delta^h$). If the distribution
$p\mapsto\Delta_p^h$ is differentiable of class $C^m$,
$m\in\field[N]\cup\{0\}$, it is said that the \emph{connection $\Delta^h$ is
(differentiable) of class} $C^m$. A connection $\Delta^h$ is of class $C^m$ if
and only if, for every $C^m$ vector field $Y$ on $E$, the vertical $Y^v$ and
horizontal $Y^h$ vector fields are of class $C^m$.

	A $C^1$ \emph{path} $\beta\colon J\to E$ is called \emph{horizontal}
(\emph{vertical}) if its tangent vector $\dot{\beta}$ is horizontal
(vertical) vector along $\beta$, \ie $\dot{\beta}(s)\in\Delta_{\beta(s)}^h$
($\dot{\beta}(s)\in\Delta_{\beta(s)}^v$) for all $s\in J$.
    A \emph{lift} $\bar\gamma\colon J\to E$ of a path $\gamma\colon J\to M$,
\ie $\pi\circ\bar\gamma=\gamma$, is called \emph{horizontal} if $\bar\gamma$
is a horizontal path, \ie when the vector field $\dot{\bar{\gamma}}$ tangent
to $\bar\gamma$ is horizontal or, equivalently, if $\dot{\bar{\gamma}}$ is
a horizontal lift of $\dot\gamma$.
    Since $\pi^{-1}(\gamma(J))$ is an $(r+1)$ dimensional submanifold of
$E$, the distribution
 $p\mapsto \Delta_p^h\cap T_p(\pi^{-1}(\gamma(J)))$ is one\ndash dimensional
and, consequently, is integrable. The integral paths of that distribution are
horizontal lifts of $\gamma$ and,
\emph{for each $p\in\pi^{-1}(\gamma(J))$, there
is a unique horizontal lift $\bar{\gamma}_p$ of $\gamma$ passing through}
$p$.

    \begin{Defn}    \label{Defn3.2}
Let $\gamma\colon[\sigma,\tau]\to M$, with $\sigma,\tau\in\field[R]$ and
$\sigma\le \tau$, and
$\bar\gamma_p$ be the unique horizontal lift of $\gamma$ in $E$ passing
through $p\in\pi^{-1}(\gamma([\sigma,\tau]))$.
The \emph{parallel transport (translation, displacement)} generated by
(assigned to, defined by) a connection $\Delta^h$ is a mapping
$\Psf\colon\gamma\mapsto\Psf^\gamma$, assigning a mapping to the path $\gamma$
    \begin{equation}    \label{3.9-4}
\Psf^\gamma\colon \pi^{-1}(\gamma(\sigma)) \to \pi^{-1}(\gamma(\tau))
\qquad \gamma\colon[\sigma,\tau]\to M
    \end{equation}
such that, for each $p\in\pi^{-1}(\gamma(\sigma))$,
    \begin{equation}    \label{3.9-5}
\Psf^\gamma(p) := \bar\gamma_p(\tau).
    \end{equation}
    \end{Defn}

    Recall now the basic properties of the parallel transports generated
by connections.

    \begin{Prop}    \label{Prop8.2}
Let
    \begin{equation}    \label{8.10-1}
\Psf\colon\gamma \mapsto \Psf^\gamma
\colon\pi^{-1}(\gamma(\sigma))\to \pi^{-1}(\gamma(\tau))
\qquad \gamma\colon[\sigma,\tau]\to M
    \end{equation}
be the parallel transport generated by a connection on some bundle
$(E,\pi,M)$. The mapping $\Psf$ has the following properties:%
\\\indent\textbf{\textup{(i)}}
    The parallel transport $\Psf$ is invariant under orientation
preserving changes of the paths' parameters. Precisely, if
$\gamma\colon[\sigma,\tau]\to M$ and
$\chi\colon[\sigma',\tau']\to[\sigma,\tau]$ is an orientation preserving
 $C^1$ diffeomorphism, then
    \begin{equation}    \label{8.11}
\Psf^{\gamma\circ\chi} = \Psf^{\gamma} .
    \end{equation}
\indent\textbf{\textup{(ii)}}
    If $\gamma\colon[0,1]\to M$ and $\gamma_{\_}\colon[0,1]\to M$ is its
canonical inverse, $\gamma_{\_}(t)=\gamma(1-t)$ for $t\in[0,1]$, then
    \begin{equation}    \label{8.12}
\Psf^{\gamma_{\_}} = \bigl(\Psf^\gamma\bigr)^{-1}.
    \end{equation}
\indent\textbf{\textup{(iii)}}
    If $\gamma_1,\gamma_2\colon[0,1]\to M$, $\gamma_1(1)=\gamma_2(0)$, and
$\gamma_1\gamma_2\colon[0,1]\to M$ is their canonical product, then
    \begin{equation}    \label{8.13}
\Psf^{\gamma_1\gamma_2} = \Psf^{\gamma_2} \circ \Psf^{\gamma_1} .
    \end{equation}
\indent\textbf{\textup{(iv)}}
    If $\gamma_{r,x}\colon\{r\}=[r,r]\to\{x\}$ for some given $r\in\field[R]$
and $x\in M$, then
    \begin{equation}    \label{8.14}
\Psf^{\gamma_{r,x}} = \id_{\pi^{-1}(x)} .
    \end{equation}
\indent\textbf{\textup{(v)}}
If
\(
\bar\gamma_p \colon s\in[\sigma,\tau]\to \bar\gamma(s)
:=
\Psf^{\gamma|\sigma,s]} (p)
\)
is the lifting of $\gamma$ through $p\in\pi^{-1}(\gamma(\sigma))$ defined by
$\Psf$ and the $C^1$ paths $\gamma_i \colon [\sigma_i,\tau_i]\to M$ are such
that
 $\gamma_1(\sigma_1)=\gamma_2(\sigma_2)$ and
 $\dot\gamma_1(\sigma_1)=\dot\gamma_2(\sigma_2)$,
then
(a) the lifted path
    \begin{equation}    \label{8.15-0}
\bar\gamma_p  \text{\textup{ is of class $C^1$}}
    \end{equation}
for any $C^1$ path $\gamma$ and
(b) the lifted paths $\bar\gamma_{1;p}$ and $\bar\gamma_{2;p}$ have equal
tangent vectors at $p$,
    \begin{equation}    \label{8.15-1}
\dot{\bar\gamma}_{1;p}(\sigma_1) = \dot{\bar\gamma}_{2;p}(\sigma_2) .
    \end{equation}
\indent\textbf{\textup{(vi)}}
If $\gamma_i \colon [\sigma_i,\tau_i]\to M$, $i=1,2$, are two $C^1$ paths and
$\gamma_1(\sigma_1)=\gamma_2(\sigma_2)$, then for every $a_1,a_2\in\field$ there
exists a $C^1$ path $\gamma_3 \colon [\sigma_3,\tau_3]\to M$ such that
$\gamma_3(\sigma_3)=\gamma_1(\sigma_1)$ ($=\gamma_2(\sigma_2)$) and the vector
tangent to the lifted path
\(
\bar\gamma_{3;p} \colon t_3\in[\sigma_3,\tau_3] \to
\bar\gamma_{3;p}(t_3) := \Psf^{\gamma_3|[\sigma_3.t_3]} (p),
\)
with $p\in\pi^{-1}(\gamma_3(\sigma_3))$, at $\sigma_3$ is
    \begin{equation}    \label{8.15-2}
\dot{\bar\gamma}_{3;p}(\sigma_3)
=
a_1 \dot{\bar\gamma}_{1;p}(\sigma_2) +  a_2 \dot{\bar\gamma}_{2;p}(\sigma_2) ,
    \end{equation}
where
\(
\bar\gamma_{i;p} \colon t_i\in[\sigma_i,\tau_i] \to
\bar\gamma_{i;p}(t_i) := \Psf^{\gamma_i|[\sigma_i.t_i]} (p) .
\)
    \end{Prop}

    \begin{Rem} \label{Rem8.1-1}
As a result of~\eref{8.11}, some properties of the parallel
transports generated by connections, like~\eref{8.12}
and~\eref{8.13}, are sufficient to be formulated/proved only for
canonical paths $[0,1]\to M$.
    \end{Rem}

    \begin{Proof}
The proofs of~\eref{8.11}--\eref{8.14} can be found in a number of works, for
example
in~\cite{K&N-1,Nash&Sen,Nikolov,Lumiste-1971,Durhuus&Leinaas,
Khudaverdian&Schwarz,Lumiste-1964,DNF-1,DNF-2,DNF-3}
    \end{Proof}

    \begin{Defn}    \label{Defn8.3}
A mapping~\eref{8.10-1} satisfying~\eref{8.11}--\eref{8.15-2} will
be called \emph{(axiomatically defined) parallel transport}.
    \end{Defn}

    \begin{Prop}    \label{Prop8.3}
Let $\Psf$ be the parallel transport assigned to a $C^m$, with
$m\in\field[N]\cup\{0\}$, connection on a \emph{smooth}, of class $C^{m+1}$,
bundle $(E,\pi,M)$. Then $\Psf$  is smooth, of class $C^{m}$, in a
sense that, if $\gamma\colon[\sigma,\tau]\to M$ is a $C^1$ path,
then $\Psf^\gamma$ is in the set of $C^m$ diffeomorphisms between
the fibres $\pi^{-1}(\gamma(\sigma))$ and
$\pi^{-1}(\gamma(\tau))$,
    \begin{equation}    \label{8.15}
\Psf\colon\gamma\mapsto\Psf^\gamma
\in \Diff^m \bigl( \pi^{-1}(\gamma(\sigma)), \pi^{-1}(\gamma(\tau)) \bigr)
\qquad \gamma\colon[\sigma,\tau]\to M .
    \end{equation}
    \end{Prop}

    \begin{Proof}
See~\cite{Steenrod,Husemoller,Greub&et_al.-1}.
    \end{Proof}

	The axiomatic approach to parallel transport was developed mainly on
the ground on the properties~\eref{8.11}--\eref{8.15} of the parallel
transports assigned to connections. However, this topic is out of the range of
the present paper and the reader is referred to the literature cited at the
beginning of section~\ref{Sect4}.


\section{Transports along paths in fibre bundles (review)}
\label{Sect4}

    The widespread approach to the concept of a ``parallel transport'' is it to
be consedered as a secondary one and defined on the basis of the
connection theory%
~\cite{K&N-1,Sachs&Wu,Nash&Sen,Warner,Bishop&Crittenden,Yano&Kon,
Steenrod,Sulanke&Wintgen,Bruhat,Husemoller,Mishchenko,R_Hermann-1,
Greub&et_al.-1,Atiyah,Dandoloff&Zakrzewski,Tamura,Hicks,Sternberg}.
However, the opposite approach, in which the parallel transport is
axiomatically defined and from it the connection theory is constructed, is
also known%
~\cite{Lumiste-1964,Lumiste-1966,Teleman,Dombrowski,Lumiste-1971,
Mathenedia-4,Durhuus&Leinaas, Khudaverdian&Schwarz,Ostianu&et_al.,Nikolov,
Poor}
and goes back to 1949~%
\footnote{~%
It seems that the earliest written accounts on this approach are
the ones due to \"U.~G.~Lumiste~\cite[sec.~2.2]{Lumiste-1964} and
C.~Teleman~\cite[chapter~IV, sec.~B.3]{Teleman} (both published
in~1964), the next essential steps being made by
P.~Dombrowski~\cite[\S~1]{Dombrowski} and W.~Poor\cite{Poor}. Besides,
the author of~\cite{Dombrowski} states that his paper is based on
unpublished lectures of prof.~Willi~Rinow~(1907--1979) in~1949;
see also~\cite[p.~46]{Poor} where the author claims that the first
axiomatical definition of a parallel transport in the tangent
bundle case is given by prof.~W.~Rinow in his lectures at the
Humboldt University in~1949. Some heuristic comments on the
axiomatic approach to parallel transport theory can be found
in~\cite[sec.~2.1]{Gromoll&et_al.} too.%
}%
; \eg it is systematically realized in~\cite{Poor}, where the
connection theory on vector bundles is investigated.
In~\cite{bp-TP-general} the concept of a ``parallel transport''
was generalize to the one of ``transport along paths''. The
relations between both concepts were analyzed
in~\cite{bp-TP-parallelT}; in particular, theorem~3.1
of~\cite[p.~13]{bp-TP-parallelT} (see theorem~\ref{Thm8.1} below) contains a
necessary and sufficient condition for a transport along paths to be
(axiomatically defined) parallel transport.

    \begin{Defn}    \label{Defn8.1}
    A \emph{transport along paths} in a (topological) bundle $(E,\pi,B)$ is a
mapping $I$ assigning to every path $\gamma\colon J\to M$ a mapping
$I^\gamma$, termed \emph{transport along} $\gamma$, such that
$I^\gamma\colon (s,t)\mapsto I^\gamma_{s\to t}$ where the mapping
    \begin{equation}    \label{8.1}
I^\gamma_{s\to t} \colon  \pi^{-1}(\gamma(s)) \to
\pi^{-1}(\gamma(t))
    \qquad s,t\in J,
    \end{equation}
called \emph{transport along $\gamma$ from $s$ to} $t$, has the properties:
    \begin{alignat}{2}  \label{8.2}
I^\gamma_{s\to t}\circ I^\gamma_{r\to s} &=
            I^\gamma_{r\to t} &\qquad  r,s,t&\in J
\\          \label{8.3}
I^\gamma_{s\to s} &= \id_{\pi^{-1}(\gamma(s))} & s&\in J ,
    \end{alignat}
where  $\circ$ denotes composition of mappings and $\id_X$ is the
identity mapping of a set $X$.
    \end{Defn}

    An analysis and various comments on this definition can be found
in~\cite{bp-TP-general,bp-TP-parallelT,bp-NF-LTP,bp-NF-D+EP}.

    As we shall see below, an important role is played by transports along
paths satisfying some additional conditions, in particular
    \begin{alignat}{2}   \label{8.5}
I_{s\to t}^{\gamma|J'} & = I_{s\to t}^{\gamma} &\qquad &s,t\in J'
\\          \label{8.6}
I_{s\to t}^{\gamma\circ \chi} & = I_{\chi(s)\to \chi(t)}^{\gamma}
&& s,t\in J^{\prime\prime},
    \end{alignat}
where $J'\subseteq J$ is a subinterval, $\gamma|J'$ is the restriction of
$\gamma$ to $J'$, and $\chi\colon J^{\prime\prime}\to J$ is a bijection of a
real interval $J^{\prime\prime}$ onto $J$.

    Putting $r=t$ in~\eref{8.2} and using~\eref{8.3}, we see that the
mappings~\eref{8.1}  are invertible and
    \begin{equation}    \label{8.4}
(I_{s\to t}^{\gamma})^{-1} = I_{t\to s}^{\gamma} .
    \end{equation}

	For discussion and more details on the transports along paths, the
reader is referred to~\cite{bp-TP-general,bp-TP-parallelT}; in particular,
there can be found the general form of these mappings, possible restrictions
on them, their relations with other differential\ndash geometric structures,
etc.


\section{A connection generated by transport along paths}
\label{Sect5}

    The following result describes how a transport along paths generates
a connection.~%
\footnote{~%
The author thanks Petko Nikolov (Physical department of the Sofia University
``St.\ Kliment Okhridski'') for a discussion which led to a correct formulation
of theorem~\ref{Thm8.01}.%
}

    \begin{Thm}    \label{Thm8.01}
Let $I$ be a transport along paths in a bundle $(E,\pi,M)$. Let
$\gamma\colon J\to M$ be a path and, for any $s_0\in J$ and
$p\in\pi^{-1}(\gamma(s_0))$,
the lift $\bar{\gamma}_{s_0,p}\colon J\to E$ of $\gamma$ be defined by
    \begin{equation}    \label{8.7}
\bar{\gamma}_{s_0,p}(t) = I_{s_0\to t}^{\gamma}(p) \qquad t\in J .
    \end{equation}
Suppose the transport $I$ is such that:\\
\indent
(a)
($C^1$ smoothness) The path
    \begin{subequations}    \label{8.70}
    \begin{equation}    \label{8.70a}
\text{$\bar{\gamma}_{s_0,p}$ \textup{is of class} $C^1$}
    \end{equation}
for every $C^1$ path $\gamma$, $s_0$ and $p$.
\\
\indent
(b) (Initial uniqueness)
If $\gamma_i \colon J_i\to M$, $i=1,2$, are two $C^1$ paths and for some
$s_i\in J_i$ is fulfilled $\gamma_1(s_1)=\gamma_2(s_2)$ and
$\dot\gamma_1(s_1)=\dot\gamma_2(s_2)$,
then the lifted paths $\bar{\gamma_i}_{;s_i,p}$, defined via~\eref{8.7}
with $p\in\pi^{-1}(\gamma_1(s_1))=\pi^{-1}(\gamma_2(s_2))$, have equal tangent
vectors at $p$,
    \begin{equation}    \label{8.70b}
\dot{\bar{\gamma}}_{1;s_1,p}(s_1)=\dot{\bar{\gamma}}_{2;s_2,p}(s_2).
    \end{equation}
%
\indent
(c) (Linearization)
If $\gamma_i \colon J_i\to M$, $i=1,2$, are two $C^1$ paths and for some
$s_i\in J_i$ is fulfilled $\gamma_1(s_1)=\gamma_2(s_2)$, then for every
$a_1,a_2\in\field$ there exists a $C^1$ path $\gamma_3 \colon J_3\to M$
(generally depending on $a_1$, $a_2$, $\gamma_1$ and $\gamma_2$) such that
$\gamma_3(s_3)=\gamma_1(s_1)$ ($=\gamma_2(s_2)$) for some $s_3\in J_3$ and the
vector tangent to the lifted path $\bar{\gamma}_{3;s_3,p}$, defined
via~\eref{8.7} with $p\in\pi^{-1}(\gamma_3(s_3))$, at $s_3$ is
    \begin{equation}    \label{8.70c}
\dot{\bar{\gamma}}_{3;s_3,p}(s_3)
=
a_1 \dot{\bar{\gamma}}_{1;s_1,p}(s_1) + a_2 \dot{\bar{\gamma}}_{2;s_2,p}(s_2).
    \end{equation}
    \end{subequations}

	Then
    \begin{multline}    \label{8.8}
\Delta^I\colon p\mapsto\Delta_p^I
:=\Bigl\{
    \frac{\od}{\od t}\Big|_{t=s_0}\bigl( \bar{\gamma}_{s_0,p}(t) \bigr)
\\
: \gamma\colon J\to M \text{\upshape\ is $C^1$ and injective},\ s_0\in J,\
    \gamma(s_0)=\pi(p)
  \Bigr\} \subseteq T_p(E) ,
    \end{multline}
with $p\in E$, is a distribution which is a connection on $(E,\pi,M)$, \ie
    \begin{equation}    \label{8.9}
\Delta_p^v\oplus\Delta_p^I = T_p(E) \qquad p\in E,
    \end{equation}
with $\Delta^v$ being the vertical distribution on $E$,
$\Delta_p^v=T_p(\pi^{-1}(\pi(p)))$.
    \end{Thm}

    \begin{Proof}
To begin with, we shall write the following result.

    \begin{Lem}    \label{Lem8.1}
If a $C^1$ path $\bar{\gamma}\colon J\to E$ is a lift of a $C^1$ path
$\gamma\colon J\to M$, $\pi\circ\bar{\gamma}=\gamma$, then
    \begin{equation}    \label{8.10}
\pi_*(\Dot{\Bar{\gamma}}(t)) = \Dot{\gamma}(t) .
    \end{equation}
    \end{Lem}

    \begin{Proof}[Proof of Lemma~\protect{\ref{Lem8.1}}]
Let $\{u^\mu=x^\mu\circ\pi,u^a\}$ be bundle coordinate system on $E$ and $p$ be
a point in its domain. Then~\eref{8.10} follows from
\(
\pi_*\bigl(\frac{\pd}{\pd u^{I}}\big|_p\bigr)
=\frac{\pd(x^\mu\circ\pi)}{\pd u^{I}}\big|_p
                    \frac{\pd}{\pd x^\mu}\big|_{\pi(p)}
\)
and
    \begin{equation}    \label{8.10-6}
\Dot{\Bar{\gamma}}^\mu
=
\Dot{\gamma}^\mu
    \end{equation}
which is a corollary of
\(
\Dot{\Bar{\gamma}}^\mu(t)
=\frac{\od(u^\mu\circ\Bar{\gamma}(t))}{\od t}
=\frac{\od(x^\mu\circ\gamma(t))}{\od t}
=\Dot{\gamma}^\mu(t)
\)
for all $t\in J$.
    \end{Proof}

	From~\eref{8.8} and~\eref{8.10}, we get
    \begin{multline*}
\pi_*(\Delta_p^I)
= \bigl\{ \Dot{\gamma}(s_0) :
     \gamma\colon J\to M \text{ is $C^1$ and injective}, s_0\in J,\
    \gamma(s_0)=\pi(p)
\bigr\}
=
T_{\pi(p)}(M)
    \end{multline*}
as $\Dot{\gamma}(s_0)$ is an arbitrary vector in
$T_{\pi(p)}(M)=T_{\gamma(s_0)}(M)$.
	Besides, the condition (c) ensures that $\Delta^I_p$ is a vector
subspace in $T_p(E)$.
    Thus $\pi_*|_{\Delta_p^I}\colon\Delta_p^I\to T_{\pi(p)}(M)$ is
a surjective mapping between vector spaces.
	It is also linear as it is a restriction of the tangent mapping $\pi_*$,
which in turn is a linear mapping~\cite[sec.~1.22]{Warner},
on the vector subspace $\Delta^I_p$.
	At last, we shall
prove that $\pi_*|_{\Delta_p^I}$ is injective, from where it follows that
$\pi_*|_{\Delta_p^I}$ is a vector space isomorphism for every
$p\in E$ which, in its turn, implies that $\Delta^I \colon p\to \Delta_p^I$ is
a distribution satisfying~\eref{8.9} as $\pi_*(\Delta_p^v)=0_{\pi(p)}\in
T_{\pi(p)}(M)$.

    If $G_i\in\Delta_p^I$, $i=1,2$, then there exist paths
$\gamma_i\colon J_i\to M$ such that $\gamma_i(s_i)=\pi(p)$ for some
$s_i\in J_i$ and $G_i=\Dot{\Bar{\gamma}}_{i;s_i,p}(s_i)$, with $i=1,2$ and the
lifted paths in the r.h.s.\ being given by~\eref{8.7}. Then
$\pi_*(G_i)=\Dot{\gamma}_i(s_i)$, due to~\eref{8.10}. Suppose that
 $\dot{\gamma}_1(s_1)=\dot{\gamma}_2(s_2)$.
Then condition~(b) entails
$\dot{\bar{\gamma}}_{1;s_1,p}(s_1)=\dot{\bar{\gamma}}_{2;s_2,p}(s_2)$
so that
	$G_1=G_2$, which means that
$\pi_*|_{\Delta_p^I}\colon\Delta_p^I\to T_{\pi(p)}(M)$ is injective.
    \end{Proof}

    \begin{Rem} \label{Rem8.0}

	The condition (a) ensures that all our constructions have a sense.
	If the condition (b) is not valid, than $\pi_*|_{\Delta^I_p}$ is linear
and surjective, but we cannot prove that it is also injective; precisely,
$\pi_*|_{\Delta^h_p}$ is injective iff the property~(b) is valid.
	At last, the condition (c) guarantees that $\Delta^I_p$ is a vector
subspace of $T_p(E)$; in fact, $\Delta^I_p$ is a vector space iff the
property~(c) is valid.
	This is quite an essential moment as the direct complement $B$ to a
vector subspace $A$ of a vector space $V$ (in our case $B=\Delta_p^h$,
$A=\Delta_p^v$ and $V=T_p(E)$) is generally \emph{not} a vector subspace; \ie
the equality $V=A\oplus B$, with $V$ and $A$ being vector spaces, does not
generally imply that $B$ is a vector space; for example, if
 $A=\{(x,0) : x\in\field[R]\}\subset\field[R]^2$, then
 $\field[R]^2=A\oplus B$ with
$B = \{(x,y)\in\field[R], y\not=0\} \cup\{(0,0)\}$,
but $B$ is not a vector subspace in $\field[R]^2$, as
 $(a,-c),(b,+c)\in B$ for $a,b,c\in\field[R]$ and $c\not=0$ but
$(a,-c)+(b,+c)=(a+b,0)\not\in B$ for $a+b\not=0$.
     \end{Rem}

    \begin{Rem}		\label{Rem8.0-1}
Since~\eref{8.70c} and~\eref{8.10} entail
    \begin{subequations}    \label{8.00}
    \begin{equation}    \label{8.00a}
\dot\gamma_3(s_3) = a_1 \dot\gamma_1(s_1) + a_2 \dot\gamma_2(s_2) ,
    \end{equation}
the path $\gamma_3$ is uniquely fixed by the conditions that it passes through
the point $\pi(p)$,
    \begin{equation}    \label{8.00b}
\gamma_3(s_3) = \pi(p) .
    \end{equation}
    \end{subequations}
Therefore, conditions (c) is equivalent to the requirement the vector tangent
at $s_3$ to the lift, given via~\eref{8.7}, of the path defined
via~\eref{8.00} to be~\eref{8.70c}.
    \end{Rem}

    \begin{Rem} \label{Rem8.1}
The role of the transport $I$ along paths in theorem~\ref{Thm8.01} is on
its base to be constructed a lifting of the paths in $M$ to paths in $E$ with
appropriate properties. Namely, such a lifting should assign to a path
$\gamma\colon J\to M$ a unique path $\Bar{\gamma}_{s_0,p}\colon
J\to E$ passing through a given point $p\in\pi^{-1}(\gamma(s_0))$,
for some $s_0\in J$, and such that
$\pi\circ\Bar{\gamma}_{s_0,p}=\gamma$ and, if $q\in
\Bar{\gamma}_{s_0,p}(t_0)$ for some $t_0\in J$, then
 $\Bar{\gamma}_{t_0,q}=\Bar{\gamma}_{s_0,p}$.
    On this ground one can generalize theorem~\ref{Thm8.01} as well
as some of the next considerations and results.
    \end{Rem}

    \begin{Defn}    \label{Defn8.2}
The connection $\Delta^I$, defined in theorem~\ref{Thm8.01}, will be
called \emph{assigned to} (\emph{defined by, generated by}) the transport $I$
along paths.
    \end{Defn}

    \begin{Rem}    \label{Rem8.00}
Note that the transport $I$ in this definition must have the
properties~\eref{8.70} described in theorem~\ref{Thm8.01}.
    \end{Rem}

	It is worth recording the simple fact that, if equation~\eref{8.5}
or~\eref{8.6} holds, then the path $\bar{\gamma}_{s_0,p}$, defined
via~\eref{8.7} is such that respectively
    \begin{alignat}{2}    \label{8.10-7}
&\bar{\gamma}_{s_0,p}|J' = (\overline{\gamma|J'}) _{s_0,p}
&\qquad &s_0\in J'
\\			  \label{8.10-8}
& \bar{\gamma}_{\chi(s_0),p} \circ \chi
=
(\overline{\gamma\circ \chi})_{s_0,p}
& \qquad &s_0\in J'' .
    \end{alignat}

    \begin{Prop}    \label{Exrc8.1}
For $C^1$ paths, the equations~\eref{8.5} and~\eref{8.6} are
consequences of the suppositions~(a)--(c) from theorem~\ref{Thm8.01} but
the converse is generally not true, i.e.~\eref{8.5} and~\eref{8.6} are
necessary but generally not sufficient conditions for~\eref{8.8} to be a
connection.
    \end{Prop}

    \begin{Proof}
	To prove~\eref{8.5}, use assumption~(b) and the local existence of a
unique $C^1$ path passing through a given point and having a fixed tangent
vector at it.
	The proof of~\eref{8.6} is more complicated. For the purpose, observe
that the paths
 $\bar\beta_{s_0,p} \colon t\to I^\beta_{s_0\to t}(p)$ and
\(
\bar\delta_{s_0,p} \colon t\to I^\gamma_{\chi(s_0)\to \chi(t)}(p)
=
\bar\gamma_{\chi(s_0),p}\circ\chi(t) ,
\)
for $s_0,t\in J''$, in $E$ are liftings of the path $\beta=\gamma\circ\chi$ and
are $\Delta^I$\ndash horizontal, \ie
$\dot{\bar\beta}_{s_0,p}(t) \in \Delta^I_{\bar\beta_{s_0,p}(t)}$ and
$\dot{\bar\delta}_{s_0,p}(t) \in \Delta^I_{\bar\delta_{s_0,p}(t)}$.
(Here  use that $\dot\beta(t)=\frac{\od \chi(t)}{\od t} \dot\gamma(\chi(t))$
and similarly for $\bar\delta_{s_0,p}$.)
	Then the existence of a unique horizontal lift of any $C^1$ path (see
also proposition~\ref{Prop8.4} below) implies
$\bar\beta_{s_0,p}=\bar\delta_{s_0,p}$.
    \end{Proof}

    \begin{Prop}    \label{Prop8.4}
If $\gamma\colon J\to M$ is an injective path and $p\in\pi^{-1}(\gamma(s_0))$
for some $s_0\in J$, then there is a unique $\Delta^I$\ndash horizontal lift
of $\gamma$ (relative to $\Delta^I$) in $E$ through $p$ and it is exactly the
path $\Bar{\gamma}_{s_0,p}$ defined by~\eref{8.7}.
    \end{Prop}

    \begin{Proof}
Simply apply the definitions~\eref{8.7} and~\eref{8.8} and use the
properties~\eref{8.1} and~\eref{8.2} of the transports along paths ($t_0\in
J$):
    \begin{multline*}
\Dot{\Bar{\gamma}}_{s_0,p} (t_0)
= \frac{\od}{\od t}\Big|_{t=t_0} \bigl( I_{s_0\to t}^{\gamma}(p) \bigr)
= \frac{\od}{\od t}\Big|_{t=t_0}
  \bigl( I_{t_0\to t}^{\gamma} ( I_{s_0\to t_0}^{\gamma}(p)) \bigr)
\\
= \frac{\od}{\od t}\Big|_{t=t_0}
  \bigl( I_{t_0\to t}^{\gamma} (\bar\gamma_{s_0,p}(t_0)) \bigr)
=
\dot{\bar\gamma}_{t_0 , \bar\gamma_{s_0,p}(t_0)}(t_0)
  \in \Delta_{\bar{\gamma}_{s_0,p}(t_0)}^I .
    \end{multline*}
    \end{Proof}

    \begin{Rem} \label{Rem8.2}
If $\gamma$ is not injective and $\gamma(s_0)=\gamma(t_0)$ for
some $s_0,t_0\in J$ such that $s_0\not=t_0$, then the paths
$\Bar{\gamma}_{s_0,p},\Bar{\gamma}_{t_0,p}\colon J\to E$ need not
to coincide as $\Bar{\gamma}_{s_0,p} = I_{t_0\to s_0}^{}\circ
\Bar{\gamma}_{t_0,p}$, due to~\eref{8.7} and~\eref{8.2}. Therefore
the $\Delta^I$\ndash horizontal lift of a non\ndash injective path
through a point in $E$ lying in a fibre over a self\ndash
intersection point of the path, if any, may not be unique. Similar
is the situation for an arbitrary connection. (This range of
problems is connected with the so\ndash called holonomy groups.)
    \end{Rem}

    \begin{Prop}    \label{Prop8.5}
The parallel transport $\Isf$ generated by the connection, defined by a
transport along paths $I$, is such that
    \begin{equation}    \label{8.20}
\Isf^\gamma = I_{\sigma\to\tau}^{\gamma}
\qquad\text{for } \gamma\colon[\sigma,\tau]\to M .
    \end{equation}
    \end{Prop}

    \begin{Proof}
According to definition~\ref{Defn3.2} and~\eref{8.7}, we have:
    \begin{equation*}
    \begin{split}
\Isf &\colon\gamma\mapsto
\Isf^\gamma\colon \pi^{-1}(\gamma(\sigma))\to \pi^{-1}(\gamma(\tau))
\qquad \gamma\colon[\sigma,\tau]\to M
\\
\Isf^\gamma &\colon p\mapsto \Isf^\gamma(p) = \Bar{\gamma}_{\sigma,p}(\tau)
= I_{\sigma\to\tau}^{\gamma}(p)
\qquad p\in\pi^{-1}(\gamma(\sigma)) .
    \end{split}
    \end{equation*}
    \end{Proof}


\section
[Links between parallel transports and transports along paths]
{Links between parallel transports and\\ transports along paths}
\label{Sect6}

	Now we shall present a modified version
of~\cite[p.~13,~theorem~3.1]{bp-TP-parallelT}.

    \begin{Thm} \label{Thm8.1}
Let $I$ be a transport along paths in a bundle $(E,\pi,M)$ and
$\gamma\colon[\sigma,\tau]\to M$. If $I$ satisfies the conditions~\eref{8.70}
from theorem~\ref{Thm8.01}, then the mapping
    \begin{equation}    \label{8.16}
\Isf\colon\gamma\mapsto \Isf^\gamma
:= I_{\sigma\to\tau}^{\gamma}
\colon \pi^{-1}(\gamma(\sigma)) \to \pi^{-1}(\gamma(\tau))
\qquad \gamma\colon[\sigma,\tau]\to M
    \end{equation}
is a parallel transport, \ie it possess the
properties~\eref{8.11}--\eref{8.15-2}, with $\Isf$ for $\Psf$.
Besides, if $I$ is smooth in a sense that
    \begin{equation}    \label{8.17}
I_{s\to t}^{{\beta}}\in
\Diff^m\bigl( \pi^{-1}(\beta(s)), \pi^{-1}(\beta(t)) \bigr)
\qquad \beta\colon J\to M   \quad s,t\in J
    \end{equation}
for some $m\in\field[N]\cup\{0\}$, then the mapping~\eref{8.16}
satisfies~\eref{8.15}, with $\Isf$ for $\Psf$.

    Conversely, suppose the mapping
    \begin{equation}    \label{8.18}
\Psf\colon\gamma\mapsto \Psf^\gamma
\colon \pi^{-1}(\gamma(\sigma)) \to \pi^{-1}(\gamma(\tau))
\qquad \gamma\colon[\sigma,\tau]\to M
    \end{equation}
is a parallel transport, \ie satisfies~\eref{8.11}--\eref{8.15-2}, and define
the mapping
    \begin{equation}    \label{8.19}
P\colon\beta\mapsto P^\beta \colon (s,t)
\mapsto P_{s\to t}^{\beta}
=
\Psf^{(\beta|[\sigma,\tau])\circ \chi_{t}^{[\sigma,\tau]}} \circ
\Bigl( \Psf^{(\beta|[\sigma,\tau])\circ \chi_{s}^{[\sigma,\tau]}}
\Bigr)^{-1}
\qquad \beta\colon J\to M,
    \end{equation}
where $s,t\in J$, $\sigma,\tau\in J$ are such that $\sigma\le\tau$ and
$[\sigma,\tau]\ni s,t$,~%
\footnote{~%
In particular, one can set $\sigma=\min(s,t)$ and $\tau=\max(s,t)$ or, if $J$
is a closed interval, define $\sigma$ and $\tau$ as the end points of $J$,
\ie $J=[\sigma,\tau]$.%
}
and
 $\chi_{s}^{[\sigma,\tau]}\colon[\sigma,\tau]\to [\sigma,s]$ are for
$s>\sigma$ arbitrary orientation preserving $C^1$ diffeomorphisms (depending
on $\beta$ via the interval $[\sigma,\tau]$).
	Then the mapping~\eref{8.19} is a transport along paths in
$(E,\pi,M)$, which transport satisfies the conditions~\eref{8.70}
, with $P$ for $I$. Besides, under the same assumptions, the
condition~\eref{8.15} for $\Psf$ implies~\eref{8.17}, with $P$ for $I$, where
$P$ is given by~\eref{8.19}.
    \end{Thm}

	\begin{Rem}	\label{Rem8.3}
Instead by~\eref{8.19}, the transport $P$ along paths generated by a
parallel transport $\Psf$ can be defined equivalently as follows. For a path
$\gamma\colon[\sigma,\tau]\to M$ and $s,t\in[\sigma,\tau]$, we put
(see~\eref{8.11})
	\begin{subequations}	\label{8.19-1}
    \begin{equation}    \label{8.19-1a}
P\colon\gamma\mapsto P^\gamma \colon (s,t)
\mapsto P_{s\to t}^{\gamma}
=
\Psf^{\gamma\circ \chi_{t}^{[\sigma,\tau]}} \circ
\Bigl( \Psf^{\gamma\circ \chi_{s}^{[\sigma,\tau]}} \Bigr)^{-1}
=
\Psf^{\gamma|[\sigma,t]}\circ \bigl(\Psf^{\gamma|[\sigma,s]}\bigr)^{-1}
\qquad \gamma\colon[\sigma,\tau]\to M
    \end{equation}
Now, for an arbitrary path $\beta\colon J\to M$, with $J$ being closed or
open at one or both its ends, we set
    \begin{equation}    \label{8.19-1b}
P\colon\beta\mapsto P^\beta \colon (s,t)
\mapsto P_{s\to t}^{\beta}
=
	\begin{cases}
P_{s\to t}^{\beta|[s,t]}	&\text{for } s\le t \\
\bigl( P_{t\to s}^{\beta|[t,s]} \bigr)^{-1}	&\text{for } s\ge t
	\end{cases}
\qquad \beta\colon J\to M   \quad s,t\in J .
    \end{equation}
	\end{subequations}
It can easily be verified that~\eref{8.19-1} are tantamount to
    \begin{equation}    \label{8.19-2}
P\colon\beta\mapsto P^\beta \colon (s,t)
\mapsto P_{s\to t}^{\beta}
=
	\begin{cases}
\Psf^{\beta|[s,t]}	&\text{for } s\le t \\
\bigl( \Psf^{\beta|[t,s]} \bigr)^{-1}	&\text{for } s\ge t
	\end{cases}
\qquad \beta\colon J\to M   \quad s,t\in J .
    \end{equation}
	\end{Rem}

    \begin{Proof}
Suppose $I$ satisfies~\eref{8.70}. To begin with, we notice
that~\eref{8.70} imply~\eref{8.5} and~\eref{8.6} by virtue of
proposition~\ref{Exrc8.1}. Equation~\eref{8.11} follows from~\eref{8.16},
and~\eref{8.6}. Equation~\eref{8.12} is a consequence of~\eref{8.16},
\eref{8.6} and~\eref{8.4} for $\gamma \colon [0,1]\to M$ and
$\tau_{\_}(t)=1-t$ for $t\in[0,1]$:
\[
\Isf^{\gamma_{\_}}
= I_{0\to1}^{\gamma\circ\tau_{\_}}
= I_{\tau_{\_}(0)\to\tau_{\_}(1)}^{\gamma}
= I_{1\to0}^{\gamma}
= (I_{0\to1}^{\gamma})^{-1}
=(\Isf^\gamma)^{-1} .
\]
To prove~\eref{8.13}, we set $\gamma_{1,2} \colon [0,1]\to M$, $\tau_1(t)=2t$
for $t\in[0,1/2]$ and $\tau_2(t)=2t-1$ for $t\in[1/2,1]$ and, applying equations
~\eref{8.16}, \eref{8.5}, \eref{8.6}, \eref{8.4} and~\eref{8.2},
we find:
    \begin{multline*}
\Isf^{\gamma_1\gamma_2}
= I_{0\to1}^{\gamma_1\gamma_2}
= I_{1/2\to1}^{\gamma_1\gamma_2} \circ I_{0\to1/2}^{\gamma_1\gamma_2}
= I_{1/2\to1}^{\gamma_1\gamma_2|[1/2,1]} \circ
  I_{0\to1/2}^{\gamma_1\gamma_2|[0,1/2]}
\\
= I_{1/2\to1}^{\gamma_2\circ\tau_2} \circ I_{0\to1/2}^{\gamma_1\circ\tau_1}
= I_{\tau_2(1/2)\to\tau_2(1)}^{\gamma_2} \circ
  I_{\tau_1(0)\to\tau_1(1/2)}^{\gamma_1}
= I_{0\to1}^{\gamma_2} \circ I_{0\to1}^{\gamma_1}
= \Isf^{\gamma_2} \circ \Isf^{\gamma_1} .
    \end{multline*}
Next, equation~\eref{8.14} is a direct corollary from~\eref{8.16}
and~\eref{8.3}. At last, the properties~(v) and~(vi) are trivial
consequences from~\eref{8.16} and
respectively~\eref{8.70a}-\eref{8.70b} and~\eref{8.70c}. So, $\Isf$
is a prarallel transport, which, evidently, satisfies~\eref{8.15} with
$\Isf$ for $\Psf$, if~\eref{8.17} is valid.

	Conversely, suppose $\Psf$ is a parallel transport and $P$ is defined
via~\eref{8.19}, or, equivalently (see remark~\ref{Rem8.3}),
via~\eref{8.19-2}. The verification of~\eref{8.1}--\eref{8.3} can be
done by applying~\eref{8.18} and~\eref{8.11}--\eref{8.14}. Therefore
$P$ is a parallel transport along paths. Similarly, the
conditions~\eref{8.70} follow from the properties~(v) and~(vi) of the
parallel transport $\Psf$. At last, the validity of~\eref{8.14}, with $P$
for $I$, is a consequence from~\eref{8.15} and~\eref{8.19-2}.
    \end{Proof}

    \begin{Defn}    \label{Defn8.4-0}
A transport $I$ along paths which has the properties~\eref{8.70}
will be called \emph{parallel} transport along paths.
    \end{Defn}

	Theorem~\ref{Thm8.1} simply says that there is a bijective
correspondence between the parallel transports along paths and the parallel
transports.

    \begin{Defn}    \label{Defn8.4}
If $I$ is a parallel transport along paths, then we say that the parallel
transport~\eref{8.16} is \emph{generated by} (\emph{defined by, assigned to})
$I$.
    Respectively, if $\Psf$ is a parallel transport, then we
say that the (parallel) transport along paths~\eref{8.19} is \emph{generated
by} (\emph{defined by, assigned to}) $\Psf$.
    \end{Defn}

    \begin{Cor} \label{Cor8.1}
The parallel transport $\Isf$ generated by the connection $\Delta^I$,
assigned to a transport $I$ along paths, is a parallel transport, \ie it
satisfies~\eref{8.11}--\eref{8.15-2} with $\Isf$ for $\Psf$.
    \end{Cor}

    \begin{Proof}
This result is a particular case of proposition~\ref{Prop8.2}. An
alternative proof can be carried out by
using~\eref{8.1}--\eref{8.70}, \eref{8.20}, and the definitions of
inverse path and product of paths. The assertion is also a
consequence of~\eref{8.20} and theorem~\ref{Thm8.1}.
    \end{Proof}

    \begin{Prop}    \label{Prop8.3-1}
If $I$ is a parallel transport along paths and $\Isf$ is the assigned to it
parallel transport, then the parallel transport along paths defined by $\Isf$
coincides with $I$. Conversely, if $\Psf$ is a parallel transport and $P$ is
the parallel transport along paths generated by $\Psf$, then the parallel
transport defined by $P$ coincides with $\Psf$.
    \end{Prop}

    \begin{Proof}
The assertions are consequences from definition~\ref{Defn8.4},
theorem~\ref{Thm8.1} and remark~\ref{Rem8.3}.
    \end{Proof}


\section
[Relations between connections and (parallel) transports along paths]
{Relations between connections and\\ (parallel) transports along paths}
\label{Sect7}

    Until this point, we have studied how a transport along paths
generates a connection (theorem~\ref{Thm8.01}) and parallel
transport (proposition~\ref{Prop8.5}). Besides,
theorem~\ref{Thm8.1} establishes a bijective correspondence
between particular class of transports along paths and mappings having (some
of) the main properties of the parallel transports generated by connections.
Below we shall pay attention, in a sense, to the opposite links, starting
from a connection on a bundle.

    \begin{Prop}    \label{Prop8.6}
Let $\Psf$ be the parallel transport assigned to a connection $\Delta^h$ on a
bundle $(E,\pi,M)$. The mapping
    \begin{subequations}    \label{8.21}
    \begin{gather}  \label{8.21a}
P\colon\gamma\mapsto P^\gamma\colon(s,t)\mapsto P_{s\to t}^{\gamma}
\qquad \gamma\colon J\to M
\\\intertext{defined by}
            \label{8.21b}
P_{s\to t}^{\gamma}
=   \begin{cases}
\Psf^{\gamma|[s,t]} &\text{for } s\le t \\
\bigl(\Psf^{\gamma|[t,s]}\bigr)^{-1}    &\text{for } s\ge t
    \end{cases}
    \end{gather}
    \end{subequations}
is a transport along paths in $(E,\pi,M)$. Moreover, $P$ is parallel
transport along paths, \ie it satisfies the equations~\eref{8.70}
 from theorem~\ref{Thm8.01} with $P$ for $I$.
    \end{Prop}

    \begin{Proof}
One should check the conditions~\eref{8.1}--\eref{8.3} and~\eref{8.70} with $P$
for $I$. The relations~\eref{8.1}--\eref{8.3} follow directly
from definition~\ref{Defn3.2} of a parallel transport generated by
a connection. The rest conditions are consequences of~\eref{8.21}
and a simple, but tedious, application of the
properties~\eref{8.11}--\eref{8.13} of the parallel transports.
	Alternatively, this proposition is a consequence of the second part
of theorem~\ref{Thm8.1} and remark~\ref{Rem8.3}.
    \end{Proof}

    \begin{Rem}	\label{Rem8.4}
Applying~\eref{8.11}--\eref{8.13}, the reader can verify that
    \begin{equation}    \label{8.22}
P_{s\to t}^{\gamma} = F^{-1}(t;\gamma)\circ F(s;\gamma)
    \end{equation}
with
    \begin{equation}    \label{8.23}
F(r;\gamma)
=   \begin{cases}
\Psf^{\gamma|[r,w]} &\text{for } r\le w \\
\bigl(\Psf^{\gamma|[w,r]}\bigr)^{-1}    &\text{for } r\ge w
    \end{cases}
\qquad r=s,t
    \end{equation}
for any (arbitrarily) fixed $w\in J$. This result is a special case of the
general structure of the transports along
paths~\cite[theorem~3.1]{bp-TP-general}.
    \end{Rem}

    \begin{Defn}    \label{Defn8.5}
The parallel transport along paths, defined by a connection $\Delta^h$ on a
bundle through proposition~\ref{Prop8.6}, will be called parallel transport
along paths \emph{assigned to (defined by, generated by)} the connection
$\Delta^h$.
    \end{Defn}

    \begin{Cor} \label{Cor8.2}
Let $\Delta^I$ be the connection generated by a parallel transport $I$ along
paths according to theorem~\ref{Thm8.01}. If $\Isf$ is the parallel
transport assigned to $\Delta^I$, then the transport along paths assigned to
$\Isf$ (or $\Delta^I$), as described in proposition~\ref{Prop8.6}, coincides
with the initial transport $I$ along paths.
    \end{Cor}

    \begin{Proof}
Substitute~\eref{8.20} into~\eref{8.21}, with $I$ for $P$ and $\Isf$ for
$\Psf$.
    \end{Proof}

    \begin{Cor} \label{Cor8.3}
Let $P$ be the transport along paths assigned to a connection $\Delta^h$
(via its parallel transport $\Psf$) according to
proposition~\ref{Prop8.6}. The connection $\Delta^P$ generated by $P$, as
described in theorem~\ref{Thm8.01}, coincides with the initial connection
$\Delta^h$, $\Delta^P=\Delta^h$.
    \end{Cor}

    \begin{Proof}
On one hand, if $p\in E$, the space $\Delta_p^P$ consists of the
vectors tangent at $s_0$ to the paths $\Bar{\gamma}_{s_0,p}\colon
t\mapsto P_{s_0\to t}^{\gamma}(p)$, with $\gamma\colon J\to M$,
$s_0\in J$, and $\pi(p)=\gamma(s_0)$, due to
theorem~\ref{Thm8.01}. On another hand, $\Delta^h_p$ consists
of the vectors tangent at $s_0$ to the paths
$\Tilde{\gamma}_{s_0,p}\colon t\mapsto \Psf^{\gamma|[s_0,t]}(p)$,
by virtue of definition~\ref{Defn3.2}. Equation~\eref{8.21b} says
that both types of paths coincide,
$\Tilde{\gamma}_{s_0,p}=\Bar{\gamma}_{s_0,p}$, so that their
tangent vectors at $t=s_0$ are identical and, consequently,
$\Delta_p^P$ and $\Delta_p^h$ are equal as sets,
$\Delta_p^P=\Delta_p^h$, for all $p\in E$.
    \end{Proof}


\section{Recapitulation}
\label{Sect8}

    Roughly speaking, the above series of results says that a connection
$\Delta^h$ is equivalent to a mapping $\Psf$ (the assigned to it parallel
transport) satisfying~\eref{8.10-1}--\eref{8.15-2} or to a mapping $P$ (the
assigned to it parallel transport along paths)
satisfying~\eref{8.1}--\eref{8.3} and~\eref{8.70} (with $P$ for $I$). Besides,
the smoothness of $\Delta^h$ is equivalent to the one of $\Psf$ or $P$. Let us
summarize these results as follows.

    \begin{Thm} \label{Thm8.2}
Given a connection $\Delta^h$ on a bundle $(E,\pi,M)$, there exists a unique
parallel transport $I$ along paths in $(E,\pi,M)$ which generates $\Delta^h$
via~\eref{8.8}, \ie $\Delta^I=\Delta^h$. Besides, the parallel transport
$\Psf$ defined by $\Delta^h$ is given by~\eref{8.16}, \ie $\Psf=\Isf$.
    \end{Thm}

    \begin{Proof}
See theorem~\ref{Thm8.01}, proposition~\ref{Prop8.5} and theorem~\ref{Thm8.1}.
    \end{Proof}

    \begin{Thm} \label{Thm8.3}
Given a parallel transport $I$ along paths in a bundle $(E,\pi,M)$, then
there exists a unique connection $\Delta^h$ on $(E,\pi,M)$ such that the
parallel transport $P$ along paths assigned to $\Delta^h$ coincides with $I$,
$P=I$. Besides, the connection $\Delta^I$ generated by $I$ is identical with
$\Delta^h$, $\Delta^I=\Delta^h$.
    \end{Thm}

    \begin{Proof}
Apply theorem~\ref{Thm8.01} and corollaries~\ref{Cor8.2}
and~\ref{Cor8.3}.
    \end{Proof}

    \begin{Thm} \label{Thm8.4}
Given a parallel transport along paths in a bundle, there is a unique
(axiomatically defined) parallel transport generating it. Conversely, given a
parallel transport, there is a unique parallel transport along paths
generating it.
    \end{Thm}

    \begin{Proof}
This statement is a reformulation of theorem~\ref{Thm8.1}.
    \end{Proof}

    \begin{Thm} \label{Thm8.5}
Given a parallel transport $\Psf$, there exists a unique
connection $\Delta^h$ generating it. Besides, the parallel
transport assigned to $\Delta^h$ coincides with $\Psf$.
    \end{Thm}

    \begin{Proof}
See theorems~\ref{Thm8.4} and~\ref{Thm8.3} and definitions~\ref{Defn3.2}
and~\ref{Defn8.4}.
    \end{Proof}

    \begin{Thm} \label{Thm8.6}
Given a connection $\Delta^h$, there is a unique parallel
transport $\Psf$ such that the defined by it parallel transport $P$ along
paths generates $\Delta^h$, $\Delta^P=\Delta^h$. Besides, $\Psf$
coincides with the parallel transport assigned to $\Delta^h$.
    \end{Thm}

    \begin{Proof}
Apply theorems~\ref{Thm8.2} and~\ref{Thm8.4}.
    \end{Proof}

    We end with the main moral of the We can summarize the above consiabove theorens. In a case of a
differentiable bundle,  the concepts ``connection,'' ``(axiomatically defined)
parallel transport'' and ``parallel transport along path'' are equivalent in a
sense that there are bijective mappings between the sets of these objects;
thus, we have the commutative diagram  shown on
figure~\vref{Equivalence-Con-PT-TP}
\renewcommand{\thefigure}{\arabic{section}.\arabic{figure}}
    \begin{figure}[ht!]
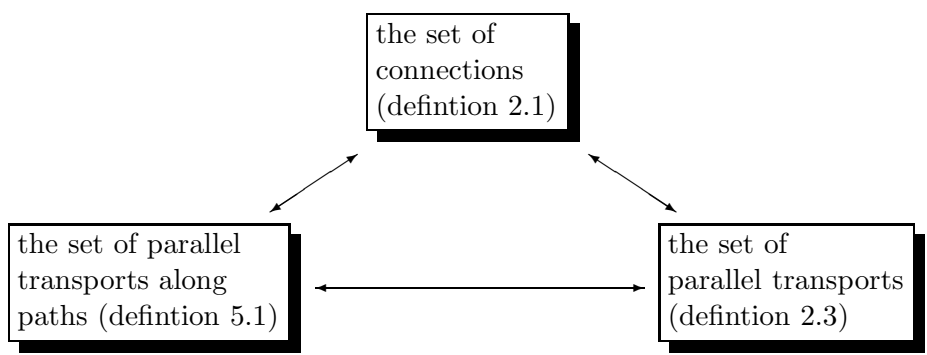

\caption
[Mappings between the sets of parallel transports, connections and parallel
transports along paths in differentiable bundles]
{\label{Equivalence-Con-PT-TP}%
Mappings between the sets of parallel transports, connections and parallel
transports along paths in differentiable bundles}
    \begin{equation*}    
    \begin{diagram}
\node[2]{\text{\shadowbox{\parbox{6.3em}{%
the set of\\ connections\\ (defintion~\ref{Defn3.1})
}}}}
\arrow{se,<>}
\arrow{sw,<>}
\\
\node{\text{\shadowbox{\parbox{9em}{
the set of parallel\\ transports along\\ paths (defintion~\ref{Defn8.4-0})
}}}}
\arrow[2]{e,<>}
\node[2]{\text{\shadowbox{\parbox{8.1em}{
the set of\\ parallel transports\\ (defintion~\ref{Defn8.3})
}}}}
    \end{diagram}
    \end{equation*}
    \end{figure}
in which the mappings indicated are described via
theorems~\ref{Thm8.2}--\ref{Thm8.6}. Besides, if one of these objects is
smooth, so are the other ones corresponding to it via the bijections
constructed in the present work.


\section{The case of topological bundles}
    \label{Subsect8.5}

	Now, following some ideas of~\cite{bp-Axiomatization-PT}, we want to say
a few words concerning bundles which are \emph{not} differentiable. The
definition~\ref{Defn3.1} of a connection is strongly related to the concept of
the tangent bundle $T(E)$ of the bundle space $E$ and, consequently, it is
senseless in a case when $E$ is a manifold of class $C^0$ of a topological
space. As we have demonstrated above (see theorems~\ref{Thm8.2}
and~\ref{Thm8.3}), a connection can equivalently be described via a suitable
parallel transport which, by definition~\ref{Defn8.3}, is a mapping possessing
the properties~\eref{8.11}--\eref{8.15-2}. One observes at first sight that the
properties~\eref{8.11}--\eref{8.14} do not depend on the differentiable
structure of the bundle and, on the contrary, the
properties~\eref{8.15-0}--\eref{8.15-2} are senseless if $E$ is not of
class $C^1$ or higher. This separation of the properties of the parallel
transports assigned to connections to ones that depend and do not depend on the
differentialble structure of the bundles points to a natural way for generalizing
the concept of a parallel transport (and hence of a connection) on bundles with
class of smoothness not higher than $C^0$.

    \begin{Defn}    \label{Defn8.6}
\index{parallel transport!in topological bundle!definition of}
\emph{A parallel transport} $\Psf$ in a \emph{topological} bundle $(E,\pi,B)$
is a mapping~\eref{8.10-1} (with $B$ for $M$) which has the
properties~\eref{8.11}--\eref{8.14}. Such a mapping will also be called
\emph{topological parallel transport} if there is a risk to mix it with the one
defined via definition~\ref{Defn8.3} in a case of a differentiable bundle.
    \end{Defn}

    \begin{Prop}    \label{Prop8.7} A parallel transport in a differentiable
bundle of class $C^1$ is a topological parallel transport.
    \end{Prop}

    \begin{Proof}
Compare definitions~\ref{Defn8.3} and~\ref{Defn8.6}.
    \end{Proof}

	Relying on our experience with connections and parallel transports in
differentiable bundles,
\emph{the concepts ``connection'' and ``(topological) transport'' in
topological bundles should be identified}
 and, hence considered as equivalent. As the following result shows, the set of
(topological) parallel transports is in a bijective correspondence with certain
subset of transports along paths in topological bundles.

    \begin{Thm} \label{Thm8.9}
Let $I$ be a transport along paths in an arbitrary topological bundle
$(E,\pi,B)$ and $\gamma\colon[\sigma,\tau]\to B$. If $I$ satisfies the
conditions~\eref{8.5} and~\eref{8.6}, then the mapping
    \begin{equation}    \label{8.16-1}
\Isf\colon\gamma\mapsto \Isf^\gamma
:= I_{\sigma\to\tau}^{\gamma}
\colon \pi^{-1}(\gamma(\sigma)) \to \pi^{-1}(\gamma(\tau))
\qquad \gamma\colon[\sigma,\tau]\to M
    \end{equation}
is a parallel transport, \ie it possess the
properties~\eref{8.11}--\eref{8.14}, with $\Isf$ for $\Psf$.

    Conversely, suppose the mapping
    \begin{equation}    \label{8.18-1}
\Psf\colon\gamma\mapsto \Psf^\gamma
\colon \pi^{-1}(\gamma(\sigma)) \to \pi^{-1}(\gamma(\tau))
\qquad \gamma\colon[\sigma,\tau]\to M
    \end{equation}
is a parallel transport, \ie satisfies~\eref{8.11}--\eref{8.14}, and define
the mapping
    \begin{equation}    \label{8.19-0}
P\colon\beta\mapsto P^\beta \colon (s,t)
\mapsto P_{s\to t}^{\beta}
=
	\begin{cases}
\Psf^{\beta|[s,t]}	&\text{for } s\le t \\
\bigl( \Psf^{\beta|[t,s]} \bigr)^{-1}	&\text{for } s\ge t
	\end{cases}
\qquad \beta\colon J\to B   \quad s,t\in J .
    \end{equation}
	Then the mapping~\eref{8.19-0} is a transport along paths in
$(E,\pi,B)$, which transport satisfies the conditions~\eref{8.5}
and~\eref{8.6}, with $P$ for $I$.
    \end{Thm}

    \begin{Proof}
To prove the first part of the theorem, we define
$\tau_- \colon t\mapsto 1-t$, $\tau_1 \colon t\mapsto 2t$ and
 $\tau_2 \colon t\mapsto 2t-1$ for $t\in[0,1]$.
Applying~\eref{8.2}--\eref{8.6} and using the notation of
proposition~\ref{Prop8.2}, we get:
    \begin{align*}
\Isf^{\gamma\circ\chi}
&= I^{\gamma\circ\chi}_{\sigma\to\tau}
 = I^{\gamma}_{\chi(\sigma)\to\chi(\tau)}
 = \Isf^{\gamma}
\\
\Isf^{\gamma_-}
&= I^{\gamma\circ\tau_-}_{0\to1}
 = I^{\gamma}_{\tau_-(0)\to\tau_-(1)}
 = I^{\gamma}_{1\to0}
 = (I^{\gamma}_{0\to1})^{-1}
 = (\Isf^{\gamma})^{-1}
\\
\Isf^{\gamma_1\gamma_2}
&= I^{\gamma_1\gamma_2}_{0\to1}
 = I^{(\gamma_1\gamma_2)|[1/2,1]}_{1/2\to1} \circ
   I^{(\gamma_1\gamma_2)|[0,1/2]}_{0\to1/2}
 = I^{\gamma_2\circ\tau_2}_{1/2\to1} \circ
   I^{\gamma_1\circ\tau_1}_{0\to1/2}
 = I^{\gamma_2}_{0\to1} \circ I^{\gamma_1}_{0\to1}
 = \Isf^{\gamma_2}\circ \Isf^{\gamma_1}
\\
\Isf^{\gamma_{r,x}}
&= I^{\gamma_{r,x}}_{r\to r} = \id_{\pi^{-1}(x)}.
    \end{align*}

	To prove the second part, we first note that~\eref{8.3} and~\eref{8.6}
follow from~\eref{8.14} and~\eref{8.11}, respectively, due to~\eref{8.19-0}.
As a result of~\eref{8.19-0}, we have, \eg for $s\le t$,
\[
P_{s\to t}^{\gamma|J'}
= \Psf^{(\gamma|J')|[s,t]}
= \Psf^{\gamma|[s,t]}
= P_{s\to t}^{\gamma}
\]
by virtue of $[s,t]\subseteq J'$. At last, we shall prove~\eref{8.2} for
$r\le s\le t$; the other cases can be proved similarly. Let us fix some
orientation-preserving bijections
 $\tau' \colon [0,1]\to [r,s]$ and
 $\tau^{\prime\prime} \colon [0,1]\to [s,t]$. Then:
    \begin{multline*}
P_{s\to t}^{\beta} \circ P_{r\to s}^{\beta}
= \Psf^{\beta|[s,t]}\circ \Psf^{\beta|[r,s]}
= \Psf^{\beta\circ\tau^{\prime\prime}}\circ \Psf^{\beta\circ\tau^{\prime}}
= \Psf^{(\beta\circ\tau^{\prime})(\beta\circ\tau^{\prime\prime})}
\\
= \Psf^{\beta\circ\tau}
= \Psf^{\beta|[r,t]}
=P_{r\to t}^{\beta}
    \end{multline*}
where
$\tau \colon a\mapsto \tau^{\prime}(2a)$  for $a\in[0,1/2]$ and
$\tau \colon a\mapsto \tau^{\prime\prime}(2a-1)$ for $a\in[1/2,1]$,
so that $\tau(0)=\tau'(0)=r$ and $\tau(1)=\tau^{\prime\prime}(1)=t$.
    \end{Proof}

    \begin{Rem}    \label{Rem8.5}
In fact, this proof is contained implicitly in the proof of
theorem~\ref{Thm8.1}.
    \end{Rem}

    \begin{Defn}    \label{Defn8.7}
\index{parallel transport along paths!in topological bundle}
In a topological bundle, a transport $I$ along paths which has the
properties~\eref{8.5} and~\eref{8.6} will be called \emph{parallel}
transport along paths. If $I$ is a parallel transport along paths, then we say
that the parallel transport~\eref{8.16-1} is \emph{generated by}
(\emph{defined by, assigned to}) $I$. Respectively, if $\Psf$ is a parallel
transport, then we say that the (parallel) transport along
paths~\eref{8.19-0} is \emph{generated by (defined by, assigned to)} $\Psf$.
    \end{Defn}

	Theorem~\ref{Thm8.9} simply says that there is a bijective
correspondence between the set of parallel transports along paths and the one
of parallel transports.

	Theorem~\ref{Thm8.9} and definition~\ref{Defn8.7} are analogues of
respectively Theorem~\ref{Thm8.1} and definition~\ref{Defn8.4-0} in a case
of a topological bundle and contain them as special cases when differentiable
bundles are considered. The last assertion is confirmed by the result of
exercise~\ref{Exrc8.1} that equations~\eref{8.5} and~\eref{8.6} are
consequences from the equations~\eref{8.70}. Thus, in terms of transports
along paths, equations~\eref{8.5} and~\eref{8.6} single out the parallel
transports which are independent of the differentiable structure of the
bundles in which they act. The same result is expressed by the
properties~\eref{8.11}--\eref{8.14} of the parallel transports assigned to
connections while~\eref{8.15-0}--\eref{8.15-2} are sensible only in
differentiable bundles and single out from the topological parallel transports
the ones defined by connections.

	We can summarize the above considerations in the commutative diagram
presented on figure~\vref{Equivalence-PT-TP}
the mappings in which are described via theorem~\ref{Thm8.9}.
\renewcommand{\thefigure}{\arabic{section}.\arabic{figure}}
    \begin{figure}[ht!]
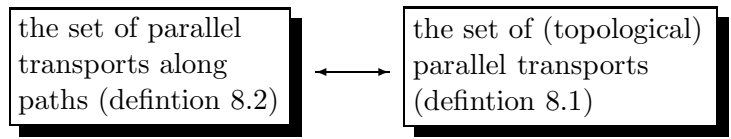

\caption
[Mappings between the sets of (topological) parallel transports and parallel
transports along paths in topological bundles]
{\label{Equivalence-PT-TP}%
Mappings between the sets of (topological) parallel transports and parallel
transports along paths in topological bundles}
    \begin{equation*}    
    \begin{diagram}
\node{\text{\shadowbox{\parbox{9em}{
the set of parallel\\ transports along\\ paths (defintion~\ref{Defn8.7})
}}}}
\arrow[2]{e,<>}
\node[2]{\text{\shadowbox{\parbox{10.1em}{
the set of (topological)\\ parallel transports\\ (defintion~\ref{Defn8.6})
}}}}
    \end{diagram}
    \end{equation*}
    \end{figure}

	As a conclusion, we can say that in topological bundles the concept of
connection does not survive and there is a bijection between the sets of
``parallel transports'' and ``parallel transports along paths'' and, in
that sense, these concepts are equivalent.


\section{The case of vector bundles}
\label{Subsect8.6}

	The purpose of this section is to be presented some links between
transports and connections in vector bundles~%
\footnote{~%
Details on the vector bundle theory can be found, for instance, in~\cite{Poor}.%
}
in a case these objects are linear,  \ie when they are compatible with the
vector structure of the bundle.

    \begin{Defn}    \label{5-Defn4.1}
A \emph{linear} transport $L$ along paths in $\field=\field[R],\field[C]$
vector bundle $(E,\pi,M)$ is a transport $L$ along paths such that
    \begin{equation}    \label{9.1}
L^\gamma_{s\to t} (\lambda u + \mu v)
=
\lambda L^\gamma_{s\to t} (u) + \mu L^\gamma_{s\to t} (v)
\qquad
\lambda,\mu\in\field, \quad u,v\in\pi^{-1}(\gamma(s))
    \end{equation}
for all paths $\gamma \colon J\to M$ and $s,t\in J$.
    \end{Defn}

    \begin{Defn} [\normalfont
			see~\protect{\cite[p.~354]{bp-NF-book}}]
		 \label{Defn9.2}
A connection on a vector bundle is called \emph{linear}
if the assigned to it parallel transport $\Psf$ is a linear mapping along every
path in the base space, \ie for all paths $\gamma\colon[\sigma,\tau]\to M$, we
have
    \begin{subequations}    \label{9.2}
    \begin{align}    \label{9.2a}
\Psf^\gamma(\rho X) &= \rho\Psf^\gamma(X)
\\		     \label{9.2b}
\Psf^\gamma(X+Y) &= \Psf^\gamma(X) + \Psf^\gamma(Y) ,
    \end{align}
    \end{subequations}
where $\rho\in\field$ and $X,Y\in\pi^{-1}(\gamma(\sigma))$.
    \end{Defn}

	Before proceeding with the topic of this section, we shall present some
results concening linear transports and linear connections.

	If $\{e_a\}$, $a,b=1,\dots,\dim\pi^{-1}(x)$, $x\in M$, is a frame in
$E$, \ie $e_a \colon M\ni x\to e_a(x)$ with $\{e_a(x)\}$ being a basis in
$\pi^{-1}(x)$, then the matrix $\Mat{L}=[L_b^a]$ of $L$ is defined by
$ L^\gamma_{s\to t} (e_a(\gamma(s))) = L_a^b(t,s;\gamma) e_b(\gamma(t)) $
and its general form is~\cite[chapter~IV, proposition~3.4]{bp-NF-book}
    \begin{equation}    \label{9.3}
\Mat{L}(t,s;\gamma) = F^{-1}(t;\gamma) F(s;\gamma)
    \end{equation}
for some non-degenerate martix-valued function $F$. The coefficients
$\Gamma_b^a(s,\gamma)$ of a $C^1$ linear transport $L$ are defined by
    \begin{equation}    \label{9.4}
\Mat{\Gamma}(s,\gamma)
:= [\Gamma_b^a(s,\gamma)]
:= \frac{\pd \Mat{L}(s,t;\gamma)}{\pd t} \Big|_{t=s}
=F^{-1}(s;\gamma)\frac{\od F(s;\gamma)}{\od s}
    \end{equation}
and determine the derivation along paths corresponding to $L$~\cite[chapter~IV,
section~3.3]{bp-NF-book}.

	If $\{X_I\}$, $I,J=1,\dots,\dim M,\dots,\dim M + \dim\pi^{-1}(x)$, $x\in
M$, is a frame in $T(E)$ adapted to
$\bigl\{ \frac{\pd }{\pd u^I} \bigr\}$ (or to $\{u^I\}$),
with $\{u^I\}$ being vector bundle coordinates in $E$, for a connection
$\Delta^h$, then~\cite[p.~368]{bp-NF-book}
    \begin{equation}    \label{9.5}
X_\mu = \frac{\pd }{\pd u^\mu} + \Gamma_\mu^\beta \frac{\pd }{\pd u^b}
\quad
X_a = \frac{\pd }{\pd u^b}
    \end{equation}
where $\mu,\nu=1,\dots,\dim M$ and $\Gamma_\mu^b$ are the (2-index) coefficients
of $\Delta^h$. According to theorem~4.1 in~\cite[p.~354]{bp-NF-book}, the
connection $\Delta^h$ is linear iff
    \begin{equation}    \label{5-4.16}
\Gamma_\mu^a(p)
= - \Gamma_{b\mu}^{a}(\pi(p)) u^b(p)
= - \bigl( (\Gamma_{b\mu}^{a}\circ \pi )\cdot u^b\bigr) (p)
\qquad p\in E ,
    \end{equation}
for some functions $\Gamma_{b\mu}^{a}$ known as the (3-index) coefficients of
$\Delta^h$.

Now we shall present some important results.

    \begin{Prop}    \label{Cor8.3-1}
In a case of a vector bundle:\\
(a) A parallel transport $I$ along paths is linear if and only if the assigned
to it parallel transport $\Isf$ is is linear;\\
(b) A parallel transport $\Psf$ is linear if and only if the generated by it
parallel transport $P$ along paths is linear.
    \end{Prop}

    \begin{Proof}
The results follow from theorem~\ref{Thm8.1} and
proposition~\ref{Prop8.3-1}.
    \end{Proof}

    \begin{Prop}    \label{Cor8.3-2}
In a case of a vector bundle:\\
	(a) The linear connections are the only connections that can be
generated by linear parallel transports along paths;\\
	(b) The linear parallel transports along paths are the only transports
along paths that can be generated by linear connections (via the assigned to
them parallel transports).
    \end{Prop}

    \begin{Proof}
The assertions are consequences from definition~\ref{5-Defn4.1} and
propositions~\ref{Cor8.3-1} and~\ref{Prop8.6}.
    \end{Proof}

    \begin{Thm}    \label{Thm8.8}
Let $\Delta^h$ be a \emph{linear} connection on a \emph{vector} bundle
$(E,\pi,M)$, $E$ and $M$ being $C^1$ manifolds. Let
$\{u^I\}=\{x^\mu\circ\pi,u^a\}$, with $\{x^\mu\}$ being a coordinate system on
$M$, be vector bundle coordinates on $E$ generated by a frame $\{e_a\}$ in $E$,
$\{X_I\}$ be the frame adapted to $\{u^I\}$ and $\Gamma_{b\mu}^{a}$ be the
3\ndash index coefficients of $\Delta^h$ in $\{X_I\}$ (see
equation~\eref{5-4.16}). If $P$ is the parallel transport along paths assigned to
$\Delta^h$, then the coefficients of $P$ in $\{e_a\}$ along a $C^1$ path
$\gamma \colon J\to M$ are ($s\in J$)
    \begin{equation}    \label{8.41}
\Gamma_{b}^{a}(s;\gamma)
=
\Gamma_{b\mu}^{a}(\gamma(s)) \dot\gamma^\mu(s).
    \end{equation}

	Conversely, if $I$ is a $C^1$ \emph{linear} transport along paths in
$(E,\pi,M)$ and in $\{e_a\}$ its coefficients have a
representation~\eref{8.41} along every $C^1$ path $\gamma$ for some
functions $\Gamma_{b\mu}^{a}$, then there exists a unique  linear connection
$\Delta^h$ on $(E,\pi,M)$ such that  the parallel transport $P$ along paths
generated by $\Delta^h$ coincides with $I$, $P=I$; besides, $\Gamma_{b\mu}^{a}$
are exactly the 3-index coefficients of $\Delta^h$.
    \end{Thm}

    \begin{Proof}
According to corollaries~\ref{Cor8.2} and~\ref{Cor8.3} and
proposition~\ref{Prop8.4}, the $\Delta^h$\ndash horizontal lift of $\gamma$
through $p\in\pi^{-1}(\gamma(s_0))$, $s_0\in J$ is given by
$\bar\gamma_{s_0,p}(s)=P_{s_0\to s}^\gamma(p)$.
	Then, in a (local) frame $\{e_a\}$ generating vector bundle coordinates
$\{u^I\}=\{u^\mu=x^\mu\circ\pi, u^a\}$, $\{x^\mu\}$ being a coordinate sysmen on
$M$, we have
\[
\bar{\gamma}_{s_0,p}(t) = P_b^a(t,s_0;\gamma) p^b e_a(\gamma(t))
\]
due to the definition of the matrix elements $P_b^a$ of $P$. Using~\eref{9.3}
and~\eref{9.4}, we find from here
$\dot{\bar{\gamma}}^\mu_{s_0,p}(t)=\dot\gamma^\mu(t)$, which is other form of
lemma~\ref{Lem8.1}, and
    \begin{equation}    \label{8.9-1}
\dot{\bar\gamma}^a_{s_0,p}(s)
=
- \Gamma_{b}^{a}(s;\gamma) \bar\gamma^b_{s_0,p}(s) .
    \end{equation}
On another hand, from the parallel transport
equation~\cite[p.~349, eq.~(3.39'a)]{bp-NF-book}, we get
\[
\dot{\bar\gamma}^a_{s_0,p}(s)
=
\Gamma_{\mu}^{a}(\bar\gamma_{s_0,p}(s))
\dot{\bar\gamma}^\mu_{s_0,p}(s),
\]
where $\Gamma_{\mu}^{a}$ are the
2-index coefficients of $\Delta^h$. Inserting here~\eref{5-4.16} and comparing
the result with the previous displayed equation, we obtain~\eref{8.41}, due
to $\pi\circ\bar\gamma_{s_0,p}=\gamma$ and~\eref{8.10}.

	Conversely, suppose~\eref{8.41} holds. Then $I$ satisfies the
conditions~(a)--(c) from theorem~\ref{Thm8.01}, so~\eref{8.8} defines a
connection $\Delta^h$. Then the parallel transport $P$ along paths generated by
$\Delta^h$ coincides with $I$, $P=I$, by virtue of corollary~\ref{Cor8.2}. As a
result of proposition~\ref{Prop8.4}, \eref{8.41} and~\eref{8.9-1}, the
$\Delta^h$ horizontal lift of $\gamma$ is such that
\[
\dot{\bar\gamma}_{s_0,p}(s)
=
- \Gamma_{b\mu}^{a}(\gamma(s))
	\dot\gamma^\mu(s) \bar\gamma_{s_0,p}^b(s) .
\]
Comparing this result with the parallel transport
equation~\cite[p.~349, eq.~(3.39'a)]{bp-NF-book}, we get (see also~\eref{8.10})
\[
\Gamma_{\mu}^{a} (\bar\gamma_{s_0,p}(s))
=
- \Gamma_{b\mu}^{a}(\gamma(s)) \bar\gamma_{s_0,p}^b(s)
\]
which entails
 $\Gamma_{\mu}^{a} = - (\Gamma_{b\mu}^{a}\circ\pi) \cdot u^b$
due to $\pi\circ\bar\gamma_{s_0,p}=\gamma$ and the arbitrariness of $\gamma$
and $p\in\pi^{-1}(\gamma(s_0))$. At last, according to
the above-cited theorem~4.1 in~\cite[p.~349]{bp-NF-book} the
connection $\Delta^h$ is linear.
    \end{Proof}

	Theorem~\ref{Thm8.8} simply means that a linear transport along
paths is a parallel transport along paths (generated by a linear connection) if
and only if equation~\eref{8.41} holds in some (and hence in any) frame
$\{e_a\}$. In this sense, a linear transport along paths is equivalent to a
linear connection iff equation~\eref{8.41} is valid in some (and hence in
any) frame $\{e_a\}$.
	We can also say that theorem~\ref{Thm8.8} describes the
subset of all linear transports along paths that can be generated by linear
connections, \viz these are the transports with coefficeints given by
equation~\eref{8.41}.

	As an application of theorem~\ref{Thm8.8}, we shall prove that the
evolution transport arising in bundle quantum mechanics~\cite{bp-BQM-1} is
generally \emph{not a parallel} linear transport along paths.~%
\footnote{~%
The below-presented proof is not quite rigorous as theorem~\ref{Thm8.8} is
proved in the finite-dimensional case while the evolution transport acts in
generally infinite dimensional Hilbert bundle.%
}
To show this, we note that the coefficients of the evolution transport
are~\cite[eq.~(2.22)]{bp-BQM-2}
    \begin{equation}    \label{9.8}
\Gamma_b^a(t;\gamma) = -\frac{1}{\iu \hbar} (\Mat{H}_\gamma^m(t))_b^a
    \end{equation}
where $\iu$ is the imaginary unit, $\hbar$ is the Planck's constant (devided
by $2\pi$) and $\Mat{H}_\gamma^m(t)$ is the matrix-bundle Hamiltonian. One can
choose a special frame which is independent of $\gamma$ and is such
that~\cite[remark~2.1]{bp-BQM-2} $\Mat{H}_\gamma^m(t)=\Mat{\mathcal{H}}(t)$,
where $\Mat{\mathcal{H}}(t)$ is the matrix of the usual Hilbert space
Hamiltonian. So, in this frame, we have
    \begin{equation}    \label{9.9}
\Gamma_b^a(t;\gamma) = -\frac{1}{\iu \hbar} \mathcal{H}_b^a(t)
    \end{equation}
with $\mathcal{H}_b^a(t)$ being the matrix elements of the system Hamiltonian
which may depend on the (time) parameter $t$ but does not depend on the path
$\gamma$. Therefore equation~\eref{8.41} is generally not valid and
consequently the evolution transport cannot be generated by a linear connection
in the general case. In particular, this conclusion is valid for all
time-independent Hamiltonians, $\frac{\pd \mathcal{H}(t)}{\pd t}=0$.


\section {Conclusion}
\label{Conclusion}

	The concept ``parallel transport'' precedes historically the one of a
``connection'' and was first clearly formulated in the
work~\cite{Levi-Civita/1917} of Levi Civita on a parallel transport of a
vector in Riemannian geometry.~%
\footnote{~%
See also the early work~\cite{Levi-Civita&Ricchi/1900}.%
}
	The connection theory was formulated approximately during the period
1920--1949 in a series of works on particular connections and their subsequent
generalizations and has obtained an almost complete form in 1950--1955 together
with the clear formulation of the concepts ``manifold''~%
\footnote{~%
The theory of manifolds was worked out by H. Whitney during the
1930s~\cite{Whitney-1936} and~\cite{de-Rham-1955} seems to be one of the first
its systematic expositions.%
}
and ``fibre bundle''~%
\footnote{~%
The first eddition of the book~\cite{Steenrod} is the first detailed and
complete presentation of fibre bundle theory. Later the book~\cite{Husemoller}
became a modern standard for references on bundle theory.%
}~\cite{Lumiste-1971}.
	It should be noted the contribution of Elie Cartan in the theory of
linear connections~\cite{Cartan-E-1924,Cartan-E-1983}.
	Connections as splitting of the second tangent
bundle $T(T(M))$ over a manifold $M$ into vertical and horizontal parts ware
first intoduced by Ehresmann in~\cite{Ehresmann-1951} and sometimes are
referred as ``Ehresmann connections''.
	The concept of connection on general fibre bundle was established at
about~1970~\cite{Liebermann-1973}.
	During that period, with a few exceptional works, the `parallel
transport' was considered as a secondary concept, defined by means of the one
of a `connection'. Later, as we pointed at the beginning of Sect.~\ref{Sect4},
there appear several attempts for axiomatizing the concept of a parallel
transport and by its means the connection theory to be constructed; e.g., this
approach is developed deeply in vector bundles in~\cite{Poor}.

	The main merit from the present paper is that now we have a necessary
and sufficient conditions when in differentiable bundles an axiomatically
defined parallel transport (or a parallel transport along paths) defines a
unique connection (with suitable properties) and \emph{vice versa}. Moreover,
the concepts ``connection,'' ``(axiomatically defined) parallel transport,''
and ``parallel transport along path'' are equivalent in a sense that there are
bijective mappings between the sets of these objects.
	In topological bundles, we have established a bijection between the
sets of parallel transports and of parallel transports along paths.
	However, the concept of a parallel transport admits a generalization to
the one of a transport along paths and there exist transports along paths that
cannot be generated by connections or parallel transports.
	In vector bundles, we have found a condition ensuring an equivalence
between linear transports and linear connections.


\section*{Acknowledgments}

	This work was partially supported by the National Science Fund of
Bulgaria under Grant No.~F~1515/2005.

	The author thanks Professor Phillip~E.\  Parker (Department of Mathematics and
Statistics, Wichita State University Wichita, Kansas, USA) for some historical
comments on connection theory.


\addcontentsline{toc}{section}{References}
\bibliography{bozhopub,bozhoref}
\bibliographystyle{unsrt}
\addcontentsline{toc}{subsubsection}{This article ends at page}

\end{document}